\documentclass[12pt]{article}
\usepackage{mathrsfs,amssymb,amsmath,amsfonts,amscd}
\def\a{\alpha}
\def\b{\beta}
\def\r{\gamma}

\def\G{\Gamma}

\def\Z{\mathbb{Z}}
\def\C{\mathbb{C}}
\def\F{\mathbb{F}}
\def\L{\mathcal{L}}
\def\W{\mathcal{W}}
\def\QED{\hfill$\Box$}

\numberwithin{equation}{section}
\newtheorem{theo}{Theorem}[section]
\newtheorem{defi}[theo]{Definition}

\newtheorem{lemm}[theo]{Lemma}

\makeatletter

\begin{document}

\title{Transposed Poisson Structures on the planar Galilean conformal algebra}
\noindent{\Large
Transposed Poisson Structures on the planar Galilean conformal algebra}\footnote{
This work was supported by National Natural Science
Foundation grants of China (11701345)}

	\bigskip
	
	 \bigskip

\begin{center}	
	{\bf
		Henan Wu\footnote{School of Mathematical Sciences, Shanxi University, Taiyuan 030006, China;
     .wuhenan@sxu.edu.cn}	  
   Wenting Zhang\footnote{ School of Mathematical Sciences, Shanxi University, Taiyuan 030006, China. 202122201017@email.sxu.edu.cn}}
\end{center}

\bigskip

\noindent{\bf Abstract:}
{\it Each $\frac{1}{2}$-derivation of the planar Galilean conformal algebra is proven to be a scalar. As a corollary, all transposed Poisson structures on the planar Galilean conformal algebra are trivial.}


\noindent{\bf MR(2010) Subject Classification:}~ 17B10, 17B65, 17B68

\bigskip

\bigskip

\section{Introduction}

The origin of Poisson algebras lies in the Poisson geometry of the 1970s. Since then, these algebras have shown their importance in several areas of mathematics and physics, such as Poisson manifolds, algebraic geometry, quantization theory, quantum groups, and classical and quantum mechanics. The study of Poisson algebras also led to other algebraic structures, such as noncommutative Poisson algebras (\cite{Casas2006Noncommutative}), Jacobi algebras, Novikov-poisson algebras (\cite{1997Novikov}). The description of all possible Poisson algebra structures with fixed Lie or associative part is an important problem in the theory of Poisson algebra (\cite{2021Poisson}).

Recently, the notion
of transposed Poisson algebra was introduced in (\cite{2020Transposed}), by exchanging the roles of the
two binary operations in the Leibniz rule defining the Poisson algebra. The relationship between $\frac{1}{2}$-derivations of Lie algebras and transposed Poisson algebras was developed in (\cite{2020}). Using this idea, possible tranposed Poisson structures with fixed Lie algebras were studied by many authors. For example,
all the possible transposed Poisson algebra structures were described on the Witt algebra and  the Virasoro algebra in (\cite{2020}), on the twisted Heisenberg-Virasoro algebra, the Schr$\ddot{o}$dinger-Virasoro algebra, the extended Schr$\ddot{o}$dinger-Virasoro algebra and the twisted Schr$\ddot{o}$dinger-Virasoro algebra in (\cite{2022}), on block Lie (super) algebras in (\cite{2022Transposed}), on Galilean and solvable Lie
algebras in (\cite{3}), on Witt type algebras in (\cite{2023Transposed}),
 on generalized Witt algebras and Block Lie algebras in (\cite{MR4617175}), respectively.

In this paper, we want to study the $\frac{1}{2}$-derivations and transposed Poisson structures of the planar Galilean conformal algebra.
This paper is organized as follows.
In Section 2, we introduce the definition of transposed Poisson algebra,
 the definition of the $\frac{1}{2}$-derivations, the definition of the planar Galilean conformal algebra.
In Section 3,  we compute the $\frac{1}{2}$-derivations of the planar Galilean conformal algebra and describe the related transposed Poisson algebra structures.

Throughout the paper,
we denote by $\C,\, \C^*,\,
\Z,\, \Z^+$ the sets of complex numbers,
nonzero complex numbers, integers, nonnegative integers, respectively. And all vector spaces and tensor products are taken over the complex field $\mathbb{C}$.

\section{Preliminaries}
In this section, we list some notations and results to be used in the paper.
The reader can refer to \cite{2020Transposed} for detail.

\begin{defi}\rm{(\cite{2020Transposed})}
let $L$ be a vector space equipped with two nonzero bilinear operations including $\cdot$ and$[,]$. The triple $(L,\cdot,[,])$ is called a transposed Poisson algebra if $(L,\cdot)$ is a commutative associative algebra and $(L, [,])$ is a Lie algebra such that for any $x,y,z\in L$,
\begin{equation}
2z\cdot[x,y]=[z\cdot x,y]+[x,z\cdot y].
\end{equation}
\end{defi}
\begin{defi}\rm{(\cite{2020})}
Let $(L,[,])$ be an algebra with multiplication $[,]$ and $\varphi$ be a linear transformation of $L$. Then $\varphi$ is called a $\frac{1}{2}$-derivation if
\begin{equation}
\varphi[x,y]=\frac{1}{2}([\varphi(x),y]+[x,\varphi(y)])
\end{equation}
for any $x,y\in L$.
\end{defi}
\begin{defi}\rm{(\cite{2016Structure})}
The planar Galilean conformal algebra $\W$ is an infinite dimensional Lie algebra over $\F$ with the basis  $\{L_m, H_m, I_m, J_m\,|\, m\in\G\}$ subject to the following nontrivial relations
\begin{equation*}
\aligned
&[L_m, L_n]=(m-n)L_{m+n},\ \ &&[L_m, H_n]=-n H_{m+n},\\
&[L_m, I_n]=(m-n)I_{m+n},\ \ &&[L_m, J_n]=(m-n)J_{m+n},\\
&[H_m, I_n]=J_{m+n},\ \ &&[H_m, J_n]=-I_{m+n}.
\endaligned
\end{equation*}
\end{defi}

\section{$\frac{1}{2}$-derivations of $\W$}
First, we develop a result on $\frac{1}{2}$-derivations of the general graded Lie algebra.
Assume $\L$ is a $G$-graded Lie algebra, i.e, $$\L=\bigoplus_{\a\in G} \L_{\a},\ \ [\L_{\a}, \L_{\b}]\subset \L_{\a+\b},$$
where $G$ is an abelian group.
Assume $D$ is an arbitrary $\frac{1}{2}$-derivation of $\L$.
For any homogenous $x_\a\in\L_{\a}$, denote $D(x_\a)=\sum y_\b$, $y_{\b}\in\L_{\b}$.
Define $D_{\r}(x_\a)=y_{\a+\r}$.
By direct computation, we have
$D_{\r}$ is a $\frac{1}{2}$-derivation of $\L$ and
$D=\sum_{\r} D_{\r}$.

Note that $\W$ is a $\Z_2\times \Z_2 \times \Z$-graded vector space with
\begin{equation*}
\W_{(\overline{0},\overline{0},m)}=\F L_m,\ \W_{(\overline{0},\overline{1},m)}=\F I_m, \W_{(\overline{1},\overline{0},m)}=\F J_m,\ \W_{(\overline{1},\overline{1},m)}=\F H_m,\ \forall \overline{0},\overline{1}\in \Z_2,\ m\in\Z.
\end{equation*}
It is obvious that $[\W_{(\overline{a_1},\overline{b_1},\r_1)}, \W_{(\overline{a_2},\overline{b_2},\r_2)}]\subset \W_{(\overline{a_1}+\overline{a_2}, \overline{b_1}+\overline{b_2},\r_1+\r_2)}$ for any $(\overline{a_1},\overline{b_1},\r_1), (\overline{a_2},\overline{b_2},\r_2)\in \Z_2\times \Z_2 \times \Z$. Denote $G=\Z_2\times \Z_2 \times \Z$. Then $\W$ is a $G$-graded Lie algebra. Now we assume $D$ is an arbitrary $\frac{1}{2}$-derivation of $\W$. Then $D=\sum D_{(\overline{a}, \overline{b},\r)}$, where $(\overline{a}, \overline{b},\r)\in G$. In the following, we attempt to characterize the homogenous $\frac{1}{2}$-derivations of $\W$ case by case. We can determine all $\frac{1}{2}$-derivations of $\W$.

Firstly we assume that $D$ is a homogenous $\frac{1}{2}$-derivation of degree $(\overline{0},\overline{0},\r)$.

\subsection{$D=D_{(\overline{0},\overline{0},\r)}$}
In this case, denote
\begin{equation*}
D(L_m)=a_mL_{m+\r},\ D(H_m)=b_mH_{m+\r},\ D(I_m)=c_mI_{m+\r},\ D(J_m)=d_mJ_{m+\r},
\end{equation*}
where $a_m,b_m,c_m,d_m\in\F$.
Since $[L_m,L_n]=(m-n)L_{m+n}$, we have
$2(m-n)D(L_{m+n})=[D(L_m),L_n]+[L_m,D(L_n)]$.
Then one gets
\begin{equation}\label{e1}
2(m-n)a_{m+n}=(m+\r-n)a_m+(m-n-\r)a_n,\ \ \forall m,n\in\Z.
\end{equation}
\begin{lemm}\label{lemm2}
The solution of (\ref{e1})
is $a_m=a_0$ for any $m\in\Z$.
\end{lemm}
{\it Proof.} Setting $n=0$ in (\ref{e1}), we have
$(m-\r)(a_m-a_0)=0$. Then $a_m=a_0$ for any $m\neq \r$.
If $\r=0$, then $a_m=a_0$ for any $m\in\Z$.
If $\r\neq0$, setting $m=\r$, $n=3\r$ in (\ref{e1}), we get $a_{\r}=a_0$.
Then the result follows.
\QED

Denote $a=a_0$.
Let $D$ act on $[L_m, I_n]=(m-n)I_{m+n}$, we have
\begin{equation}\label{e2}
2(m-n)c_{m+n}=(m+\r-n)a+(m-n-\r)c_n,\ \ \forall m,n\in\Z.
\end{equation}
Similar to the proof of Lemma \ref{lemm2}, we have
$c_m=a$ for any $m\in\Z$.
Let $D$ act on $[L_m, J_n]=(m-n)J_{m+n}$, we have
\begin{equation}\label{e3}
2(m-n)d_{m+n}=(m+\r-n)a+(m-n-\r)d_n,\ \ \forall m,n\in\Z.
\end{equation}
Similar to the proof of Lemma \ref{lemm2}, we have
$d_m=a$ for any $m\in\Z$.

Let $D$ act on $[L_m, H_n]=-n H_{m+n}$, we have
\begin{equation}\label{e4}
-2nb_{m+n}=-na+(-n-\r)b_n,\ \ \forall m,n\in\Z.
\end{equation}
Setting $n=1$, we have $-2b_{m+1}=-a+(-1-\r)b_1$, implying $b_m=b_n$ for any $m,n\in\Z$. Denote $b=b_m$. Then we get $(\r-n)b=-na$ for any $n\in\Z$.
Consequently, $a=b=0$ if $\r\neq 0$ and $a=b$ if $\r=0$.

So we can get $D(L_m)=aL_{m+\r},\ D(H_m)=aH_{m+\r},\ D(I_m)=aI_{m+\r},\ D(J_m)=aJ_{m+\r}$ when $D=D_{(\overline{0},\overline{0},\r)}$.

Then we assume that $D$ is a homogenous $\frac{1}{2}$-derivation of degree $(\overline{0},\overline{1},\r)$.

\subsection{$D=D_{(\overline{0},\overline{1},\r)}$}
In this case, denote
\begin{equation*}
D(L_m)=x_m I_{m+\r},\ D(H_m)=y_m J_{m+\r},\ D(I_m)=z_m L_{m+\r},\ D(J_m)=w_m H_{m+\r},
\end{equation*}
where $x_m, y_m, z_m, w_m\in\F$.
Let $D$ act on $[L_m, L_n]=(m-n)L_{m+n}$, we have
\begin{equation}\label{e5}
2(m-n)x_{m+n}=(m+\r-n)x_m+(m-n-\r)x_n,\ \ \forall m,n\in\Z.
\end{equation}
Similar to the proof of Lemma \ref{lemm2}.
We have $x_m=x_n$ for any $m,n\in\Z$.
Denote $x=x_m$.
Let $D$ act on $[L_m, H_n]=-n H_{m+n}$, we have
\begin{equation}\label{e6}
-2ny_{m+n}=-x+(m-n-\r)y_n,\ \ \forall m,n\in\Z,
\end{equation}
which implies $x=0, y_m=0, m\in\Z$.

Let $D$ act on $[L_m, I_n]=(m-n)I_{m+n}$, we have
\begin{equation}\label{e7}
2(m-n)z_{m+n}=(m-n-\r)z_n,\ \ \forall m,n\in\Z.
\end{equation}
\begin{lemm}\label{lemm4}
The solution of (\ref{e7})
is $z_m=0$ for any $m\in\Z$.
\end{lemm}
{\it Proof.}
Similar to the proof of Lemma \ref{lemm2},
 Setting $m=0$ in (\ref{e7}), we have
$-2nz_n=(-n-\r)z_n$. We can get $(n-\r)z_n=0$. Then $z_n=0$ for any $n\neq \r$.
If $\r=0$, then $z_n=0$ .
If $\r\neq0$, setting $m=-\r$, $n=\r$ in (\ref{e7}), we get $z_{\r}=0$.
 We have
$z_m=0$ for any $m\in\Z$.

Similar to the proof of Lemma \ref{lemm2}.
Let $D$ act on $[L_m, J_n]=(m-n)J_{m+n}$, we have
\begin{equation}\label{e8}
2(m-n)w_{m+n}=(-n-\r)w_n,\ \ \forall m,n\in\Z.
\end{equation}

\begin{lemm}\label{lemm5}
The solution of (\ref{e8})
is $w_m=0$ for any $m\in\Z$.
\end{lemm}
{\it Proof.}
Setting $m=0$ in (\ref{e8}), we have
$(n-\r)w_n=0$ Then $w_n=0$ for any $n\neq \r$.
If $\r=0$, then $w_n=0$ .
If $\r\neq0$, setting $m=-\r$, $n=\r$ in (\ref{e8}), we get $w_{\r}=0$,
$w_m=0$ for any $m\in\Z$.

 Now we can get $D(L_m)=0,\ D(H_m)=0,\ D(I_m)=0,\ D(J_m)=0$ when $D=D_{(\bar{0},\bar{1},\r)}$.

 Then we assume that $D$ is a homogenous $\frac{1}{2}$-derivation of degree $(\overline{1},\overline{0},\r)$.

\subsection{$D=D_{(\overline{1},\overline{0},\r)}$}
In this case, denote
\begin{equation*}
D(L_m)=h_mJ_{m+\r},\ D(H_m)=e_mI_{m+\r},\ D(I_m)=f_mH_{m+\r},\ D(J_m)=g_mL_{m+\r},
\end{equation*}
where $h_m,e_m,f_m,g_m\in\F$.

Let $D$ act on $[L_m, L_n]=(m-n)L_{m+n}$, we have
\begin{equation}\label{e9}
2(m-n)h_{m+n}=(m+\r-n)h_m+(m-n-\r)h_n,\ \ \forall m,n\in\Z.
\end{equation}

Similar to the proof of Lemma \ref{lemm2}, we have $h_m=h_n$ for any $m,n\in\Z$.
Denote $h=h_m$.
Let $D$ act on $[L_m, H_n]=-n H_{m+n}$, we have
\begin{equation}\label{e10}
-2ne_{m+n}=h+(m-n-\r)e_n,\ \ \forall m,n\in\Z,
\end{equation}

which implies $h=0, e_m=0, m\in\Z$.
Let $D$ act on $[L_m, I_n]=(m-n)I_{m+n}$, we have
\begin{equation}\label{e11}
2(m-n)f_{m+n}=(-n-\r)f_n,\ \ \forall m,n\in\Z.
\end{equation}

Similar to the proof of Lemma \ref{lemm5},
 Setting $m=0$ in (\ref{e11}), we have
$(n-\r)f_n=0$. Then $f_n=0$ for any $n\neq \r$.
If $\r=0$, then $f_n=0$ .
If $\r\neq0$, setting $m=-\r$, $n=\r$ in (\ref{e11}), we get $f_{\r}=0$.
 We have
$f_m=0$ for any $m\in\Z$.

Let $D$ act on $[L_m, J_n]=(m-n)J_{m+n}$, we have
\begin{equation}\label{e12}
2(m-n)g_{m+n}=(m-n-\r)g_n,\ \ \forall m,n\in\Z.
\end{equation}

Similar to the proof of Lemma \ref{lemm4}. Setting $m=0$ in (\ref{e12}), we have
$(n-\r)g_n=0$. Then $g_n=0$ for any $n\neq \r$.
If $\r=0$, then $g_n=0$.
If $\r\neq0$, setting $m=-\r$, $n=\r$ in (\ref{e12}), we get $g_{\r}=0$,
$g_m=0$ for any $m\in\Z$.

So we can get $D(L_m)=0,\ D(H_m)=0,\ D(I_m)=0,\ D(J_m)=0$ when $D=D_{(\overline{1},\overline{0},\r)}$.

Then we assume that $D$ is a homogenous $\frac{1}{2}$-derivation of degree $(\overline{1},\overline{1},\r)$.
\subsection{$D=D_{(\overline{1},\overline{1},\r)}$}
In this case, denote
\begin{equation*}
D(L_m)=i_m H_{m+\r},\ D(H_m)=j_mL_{m+\r},\ D(I_m)=k_mJ_{m+\r},\ D(J_m)=l_mI_{m+\r},
\end{equation*}
where $i_m,j_m,k_m,l_m\in\F$.
$D$ restricts to a $1/2$-derivation of the Heisenberg-Virasoro algebra,

Let $D$ act on $[L_m, H_n]=-n H_{m+n}$, we have
\begin{equation}\label{e13}
-2nj_{m+n}=(m-n-\r)j_n,\ \ \forall m,n\in\Z,
\end{equation}

\begin{lemm}\label{lemm6}
The solution of (\ref{e13})
is $j_m=0$ for any $m\in\Z$.
\end{lemm}
{\it Proof.}
Setting $m=0$ in (\ref{e13}), we have
$(n-\r)j_n=0$. Then $j_n=0$ for any $n\neq \r$.
If $\r\neq 0$, then $j_0=0$ .
If $\r\neq0$, setting $m=-\r$, $n=\r$ in (\ref{e13}), we get $j_{\r}=0$. If $\r=0$ in (\ref{e13}) we have
\begin{equation}\label{e14}
2nj_{m+n}=(n-m)j_n,\ \ \forall m,n\in\Z,
\end{equation}
Setting $m=0$ in (\ref{e14}), we get if $n\neq 0$, then $j_{n}=0$. Setting $n=-m=1$ in (\ref{e14}), we have $j_{1}=0$, so $j_{0}=0$.
 We have
$j_m=0$ for any $m\in\Z$.

Let $D$ act on $[L_m, I_n]=(m-n)I_{m+n}$, we have
\begin{equation}\label{e15}
2(m-n)k_{m+n}=i_{m}+k_{n}(m-n-\r),\ \ \forall m,n\in\Z.
\end{equation}
\begin{lemm}\label{lemm6}
The solution of (\ref{e15})
is $k_n=0$ for any $n\in\Z$.
\end{lemm}
{\it Proof.}
 Setting $m=0$ in (\ref{e15}), we have
$(\r-n)k_{n}=i_{0}$. Then $i_{0}=k_{n}=0$ for any $n\in Z$.

Let $D$ act on $[L_m, J_n]=(m-n)J_{m+n}$, we have
\begin{equation}\label{e16}
2(m-n)l_{m+n}=-i_{m}+(m-n-\r)l_{n},\ \ \forall m,n\in\Z.
\end{equation}
 Setting $m=0$ in (\ref{e16}), we have
$(n-\r)l_{n}=i_{0}$ Then $i_{0}=l_{n}=0$ for any $n\in Z$.

Let $D$ act on $[L_m, L_n]=(m-n)L_{m+n}$, we have
\begin{equation}\label{e17}
2(m-n)i_{m+n}=(m+\r)i_m+(-n-\r)i_n,\ \ \forall m,n\in\Z.
\end{equation}

\begin{lemm}\label{lemm10}
The solution of (\ref{e17})
is $i_m=0$ for any $m\in\Z$.
\end{lemm}
{\it Proof.}
Setting $m=0$ in (\ref{e17}),we have
\begin{equation}\label{e18}
-2ni_{n}=\r i_0+(-n-\r)i_n,\ \ \forall m,n\in\Z.
\end{equation}
We also have
\begin{equation}\label{e19}
i_{n}=\frac{\r}{\r-n}i_{0},\ \ n\neq \r.
\end{equation}

Setting $n=-m$, taking $m\neq \{\r,-\r,0\}$ in (\ref{e17}), we have
\begin{equation}\label{e20}
4mi_{0}=\frac{4m\r^2}{(\r-m)(\r+m)}i_{0},\ \ n\neq \r.
\end{equation}
By(\ref{e20}) we can get $i_{0}=0$ for any $m\in Z$.
If $n\neq \r$ in(\ref{e19}), we have $i_{n}=0$. setting $m=-\r$, $n=\r$ in (\ref{e17}), If $\r\neq 0$ we get $i_{\r}=0$. By (\ref{e15}) we have $(\r-n)k_{n}=i_{0}$, so when $\r=0$, we have $i_{0}=0$ for any $n\in Z$. We have $i_m=0$ for any $m\in\Z$.

So we can get $D(L_m)=0,\ D(H_m)=0,\ D(I_m)=0,\ D(J_m)=0$ when $D=D_{(\overline{1},\overline{1},\r)}$.

 As a summary, each $\frac{1}{2}$-derivation of the planar Galilean conformal algebra is proven to be a scalar. Due to Theorem 8(\cite{2020}), we get
\begin{theo}
The transposed Poisson structures on $\W$ are trivial.
\end{theo}

\bibliographystyle{unsrt}
\bibliography{planar}

\begin{thebibliography}{10}

\bibitem{Casas2006Noncommutative}
Casas, JM, and Datuashvili.
\newblock Noncommutative leibniz-poisson algebras.
\newblock {\em COMMUN ALGEBRA}, 2006.

\bibitem{1997Novikov}
Xiaoping Xu.
\newblock Novikov–poisson algebras.
\newblock {\em Journal of Algebra}, 190(2):253--279, 1997.

\bibitem{2021Poisson}
H.~Albuquerque, E.~Barreiro, S.~Benayadi, M.~Boucetta, and J.M. Sánchez.
\newblock Poisson algebras and symmetric leibniz bialgebra structures on oscillator lie algebras.
\newblock {\em North-Holland}, 2021.

\bibitem{2020Transposed}
C.~Bai, R.~Bai, L.~Guo, and Y.~Wu.
\newblock Transposed poisson algebras, novikov-poisson algebras and 3-lie algebras.
\newblock 2020.

\bibitem{2020}
Bruno Leonardo~Macedo Ferreira, Ivan Kaygorodov, and Viktor Lopatkin.
\newblock $\\frac{1}{2}$-derivations of lie algebras and transposed poisson algebras.
\newblock 2020.

\bibitem{2022}
Lamei Yuan and Qianyi Hua.
\newblock -(bi)derivations and transposed poisson algebra structures on lie algebras.
\newblock {\em Linear and Multilinear Algebra}, 70(22):7672--7701, 2022.

\bibitem{2022Transposed}
Ivan Kaygorodov and Mykola Khrypchenko.
\newblock Transposed poisson structures on block lie algebras and superalgebras.
\newblock 2022.

\bibitem{3}
Ivan Kaygorodov, Viktor Lopatkin, and Zerui Zhang.
\newblock Transposed poisson structures on galilean and solvable lie algebras.
\newblock {\em arXiv e-prints}, 2022.

\bibitem{2023Transposed}
Ivan Kaygorodov and Mykola Khrypchenko.
\newblock Transposed poisson structures on witt type algebras.
\newblock {\em Linear Algebra and its Applications}, 2023.

\bibitem{MR4617175}
Ivan Kaygorodov and Mykola Khrypchenko.
\newblock Transposed {P}oisson structures on generalized {W}itt algebras and {B}lock {L}ie algebras.
\newblock {\em Results Math.}, 78(5):Paper No. 186, 20, 2023.

\bibitem{2016Structure}
Shoulan Gao, Dong Liu, and Yufeng Pei.
\newblock Structure of the planar galilean conformal algebra.
\newblock {\em Reports on Mathematical Physics}, 78(1):107--122, 2016.

\end{thebibliography}

\end{document}